\documentclass[12pt,a4paper]{article}
\usepackage[utf8]{inputenc}
\usepackage{amsmath}
\usepackage{amsthm}
\usepackage{amsfonts}
\usepackage{amssymb}
\usepackage{url}
\usepackage{graphicx}
\graphicspath{ {images/} }

\newtheorem{remark}{Remark}
\newtheorem{lemma}{Lemma}
\newtheorem{corollary}{Corollary}
\newtheorem{theorem}{Theorem}
\newtheorem*{theorem*}{Theorem}

\newtheorem{definition}{Definition}
\author{Alexander Kushkuley \\ kushkuley@gmail.com}
\title{Some Remarks on Random Vectors and $O(n)$-Invariants}
\begin{document}
\maketitle

\begin{abstract}
\noindent  Computations involving invariant random vectors  are directly related to the theory of invariants (cf. e.g \cite{Weing_1}). 
Some  simple observations along these lines are presented in this paper. We note in particular  that sum of  elements of the standard basis of $  O(n)$-invariants is equal to the expectation of a random Veronese tensor up to a known scalar multiplier. 
\end{abstract}

\section{Introduction}
 We start with a brief overview of objects we will be dealing with (similar scenario was introduced in \cite{irr}).  Let a compact Lie group $G$ act on a real vector space $V$. We  will always  assume that standard Euclidean  scalar product $ <,> $ on $V$  is $G$-invariant (cf. e.g \cite{Adams}).
 
%

\bigskip\noindent Let $ \mathcal{U} \subset V $ be a $G$-invariant subset. A probability measure $\mu$ defined on $\mathcal{U} $ is called $G$-invariant if $ \mu (g\mathcal{U}') = \mu(\mathcal{U}') $ for any measurable  $\mathcal{U}' \subset \mathcal{U} $ and any $ g \in G $  (cf. e.g. \cite{Federer}, \cite{Leng}). $\; \mathcal{U}$ will be called an invariant probability subspace of $ V $ if such a measure exists. 
 
 \bigskip\noindent 
 A random vector $x$  distributed on $\mathcal{U} $ according to the invariant measure $\mu$ will be called $G$-invariant random vector or just \emph{invariant random vector}. 

\bigskip\noindent 
 The expectation of a random vector $x$    defined on a probability space  $ \mathcal{U} \subset V $ will be denoted by $ \mathbb{E}_{\mu}(x) $ or  just $ \mathbb{E} (x) $. It is clear from definitions that for invariant random vectors  $ \mathbb{E}(gx) =  \mathbb{E}(x)  $.

\bigskip\noindent 
 Let $ \chi $ be a normalized Haar measure on $G $ (cf. e.g. \cite{Adams}-\cite{Leng}).
An invariant projector $ P^G : V \rightarrow V^G  $ onto the subspace of $G$-fixed points in $V$ can be thus written as
\begin{equation}
P^G(v) = \int_G gv d\chi(g)  \nonumber
\end{equation}
for any $ v \in V $ (cf. e.g. \cite{Adams}).

\bigskip\noindent  
The following well known statement is  a starting point for this discussion (cf. e.g  \cite{Weing_1} or \cite{irr})

\begin{lemma}  If $x$ is a random vector defined on any invariant probability subspace of a real representation $V$ of a compact (Lie) group $G$, then
	
\begin{equation}
	\mathbb{E}(x) = \mathbb{E}( P^G (x) )  \nonumber
\end{equation} 
\end{lemma}

\subsection{Tensors and Invariants}
Let $G$ be a compact Lie group 
and $V$ be a (linear space of) real $G$-representation of dimension $n$.
Let $V^{\otimes,m}$ be a space of  $m$-tensors over $V$. 
Vector space  $V^{\otimes,m}$ has a natural $G$-invariant dot-product defined by a rule :
\begin{equation}
	< x_1 \otimes \cdots
	\otimes x_m , \; y_1 \otimes \cdots
	\otimes y_m > \; = \; <x_1, y_1> \cdots   \nonumber
	<x_m, y_m> 
\end{equation}
where $ 
 x_1 \otimes \cdots
\otimes x_m , \; y_1 \otimes \cdots
\otimes y_m $ are any two decomposable tensors in   $V^{\otimes,m}$.
We fix once and for all an orthonormal basis $ e_1, e_2, \cdots e_n $ in $V$.
\begin{definition}
Let $\mathcal{P} = \{J_1, J_2, \cdots, J_k \} $ be a partition of the set of indexes $ J= \{1,2, \cdots, m \} $.  For $ v_1, \cdots v_k \in V $ denote by 
\begin{equation}
 v_1(J_1) \otimes \cdots \otimes v_k(J_k) 
\end{equation}
a decomposable tensor in $V^{\otimes,m}$ that has vector $v_1$  at indexes enumerated by $J_1$, vector 
$ v_2$ at positions enumerated by   $J_2$ and so on.
\end{definition} 
\noindent For example, if $m = 5,
J_1 = \{1,3,5\} $ and $ J_2 = \{2,4\} $ then $ u(J_1) \otimes v(J_2) $ denotes the tensor  $ u \otimes v \otimes u \otimes v\otimes u $. 
\begin{remark}
	The expression (1) should not be understood as a tensor product. It is used for notational purposes only. More standard way to specify tensors of this kind could employ Young tableau and permutation group, but the definition given above is sufficient for our purposes. In any case, assigning a meaning to an expression (1) is beyond the scope of this paper
\end{remark}
\noindent We will always assume that all the parts of  a partition $\mathcal{P} = \{J_1, J_2, \cdots, J_k \} $  are non-trivial ($ |J_i|\neq 0, \; j = 1, \cdots k $) and we will use the notation  $k = |\mathcal{P}| $ for the number of parts of the partition $\mathcal{P}$. 
 If $ m = 2k $ and all the part sizes $|J_i|, \; i = 1, \; \cdots, k $ are equal, the partition $\mathcal{P}$ is called  \emph{pairing} (cf. e.g.  \cite{Weing_1}- \cite{gram}). 
 In case of $ m = 2$ there is just  one pairing $ J_0 = (1,2) ,  \; e_i(J_0) = e_i \otimes e_i, \; i = 1,\cdots, n $ and there is a $G$-invariant tensor 
 \begin{equation}
  I( J_0) = 	\sum_{i=1}^n e_i(J_0) = 	
 	\sum_{i=1}^n e_i \otimes e_i \; \in V\otimes V  \nonumber
 \end{equation} 
 In general, for a pairing $\mathcal{P} = \{J_1, J_2, \cdots, J_k \} $ on $J_{m=2k}$ there is a corresponding $G$-invariant tensor
\begin{equation}
 I(\mathcal{P}) = 	\sum_{i_1,i_2, \cdots , i_k \; = \;  1}^{n}  e_{i_1}(J_{1}) 
	\otimes e_{i_2}(J_{2}) \otimes \cdots e_{i_k}(J_{k})  \nonumber
\end{equation}
Let $ \Pi_m$ denotes the set of all pairings on the set of indexes $J$. For the full orthogonal group $O(V) = O(n) $  the classical First Fundamental Theorem (FFT) of  Invariant Theory  (see \cite{Lehrer} for example) actually states that   
\begin{theorem*}[\bf{FFT for Orthogonal Group}]
 There are no non-trivial $O(V) $ invariants in   $V^{\otimes,m}$ if $ m $ is odd and 
 if $m=2k$, the invariants $I(\mathcal{P}), \mathcal{P} \in \Pi_m $ span the linear space $ \mathcal{I}_m $ of $ O(V) $ invariants in $V^{\otimes,m}$.
\end{theorem*}
\noindent There are $ (2k -1) !!  \equiv  3 \cdot 5 \cdots (2k - 1 ) $ pairings on the set of size $m=2k$.  Take any pairing, e.g. $ \mathcal{P}_0 =  (1,2)(3,4) \cdots ( m-1 = 2k-1, m = 2k) $. The permutation group $S_m$ of the index set $J_m$ naturally acts on the set of  pairings  $\Pi_m$  and it is clear that   $\Pi_m$   is an $S_m$ orbit of $ \mathcal{P}_0 $. It should be also clear that 
\begin{lemma} Given two pairings $ \mathcal{P}_1,  \mathcal{P}_2 \in \Pi_m $ 
\begin{equation}
	< I(\sigma \mathcal{P}_1),  I(\mathcal{P}_2) > \; = \; < I(\mathcal{P}_1),  I( \sigma^{-1}\mathcal{P}_2) >  \nonumber
\end{equation} 
			for any $ \sigma \in S_m $.
\end{lemma}
\noindent   Let  $  S_k  \approx \frak{S}_k \subset S_{m=2k} $ be a (sub)group of permutations of even-numbered indexes.  It is easy to see, that if $ k > n $ then  
\begin{equation}
	\sum_{\sigma \in \frak{S}_k}  \text{sign}( \sigma) \; I(\sigma \mathcal{P}_0) = 0
\end{equation}
\noindent where $\text{sign}( \sigma) = \pm 1 $ is a sign of a permutation. All linear relations between $O(V)$ invariants are described by the classical Second Fundamental Theorem (SFT) (see for example a concise exposition in  \cite{Lehrer})

\begin{theorem*}[\bf{SFT for Orthogonal Group}]
The elements of the set of invariants $ \{ I(\mathcal{P}), \; \mathcal{P} \in \Pi_{m=2k} \}$ are linearly independent if $k\leq n $.
	For $ k > n $ all linear relations between elements of this set  
	 are linear combinations of the identities of the form (2) above.
\end{theorem*}
\begin{definition} Call the generating set of $O(n)$-invariants  $I(\mathcal{P}), \mathcal{P} \in \Pi_{m=2k} $ the \b{standard set} (basis if $n \geq k $) of the space $ \mathcal{I}_m $ of $ O(V) $ invariants.
 The elements of this set will be called standard invariants
\end{definition}

\subsection{Random Tensors and Veronese Surface}
For a unit sphere $ S(V) \subset V$ one has an equivariant  Veronese map 
\begin{equation}
S(V) \ni	  x \; \xrightarrow{\nu} \; x ^{\otimes,m} = x \otimes \cdots \otimes x \in 
	S(V^{\otimes,m})   \nonumber
\end{equation} 
 The image $\mathfrak{V}_m $ of  the map $\nu$ is called  \emph{Veronese} surface (cf. e.g. \cite{Veronese}).  
If $ \mathcal{U} \subset S(V)$ is a  $G$-invariant probability subspace in $S(V) $ then   $ \nu(\mathcal{U})  \subset \mathfrak{V}_m $  is  an invariant probability subspace of the Veronese surface (cf. e.g. \cite{irr}) and for an invariant  random unit vector 
$ x \in \mathcal{U}$ there is a corresponding 
invariant random Veronese tensor  $ x ^{\otimes,m} \in \mathcal{U}^{\otimes,m}  \subset \mathfrak{V}_m $.  
As in \cite{irr} we have   
\begin{lemma}
	For   independent invariant random vectors $x,y$ defined on any  invariant probability subspace  $ \mathcal{U} \subset S(V) $
		\begin{equation}
		\mathbb{E}( <x,y>^m )  \; = \;  	
		 \mathbb{E} ( \parallel P^G( 
		 x ^{\otimes,m} ) \parallel^2 )    
	\end{equation}
\end{lemma}
\noindent Proof. By Lemma 1 and independence of $x$ and $y$ 
	\begin{equation}
	\mathbb{E}( <x,y>^m ) = \mathbb{E} (<x ^{\otimes,m} ,\; y ^{\otimes,m} > )   \; = \; < \mathbb{E} ( x ^{\otimes,m} ) ,\;  \mathbb{E} ( y ^{\otimes,m} ) > \;  =  \; \mathbb{E} ( \parallel P^G( 
	x ^{\otimes,m} ) \parallel^2 )    
	 \nonumber
\end{equation}
From Lemma 1 (see  \cite{irr} for details),  we have also
\begin{lemma} 
If the representation of $G$ in $V$ is irreducible then for a random vector $x$ defined on any  invariant probability subspace  $ \mathcal{U} \subset S(V) $
\begin{equation}                                            
	\mathbb{E}( x \otimes x  ) = \frac{1}{n} I(J_0) \;  \equiv 	
\frac{1}{n}	\sum_{i=1}^n e_i \otimes e_i  \nonumber
 \end{equation}
		
\end{lemma}

\noindent Moreover, in accordance with Definition 1 (cf. (1)), for any partition $\mathcal{P} = \{J_1, J_2, \cdots, J_l \} $  of the set of indexes $ J= \{1,2, \cdots, m \} $ and for any 
invariant probability subspace  $ \mathcal{U} \subset S(V) $
there is a generalized Veronese surface defined as an image of the map 
\begin{equation}
 \mathcal{U}^{\times,l} \ni	 (x_1, x_2, \cdots, x_l)  \rightarrow
  x_1(J_1) \otimes \cdots \otimes  x_l(J_l) 
\end{equation}
and applying Lemmas 1 and  4 we get 
\begin{lemma}
For any set of independent random vectors
	$ (x_1, x_2, \cdots x_l) \in \mathcal{U} $ 
	defined on any  invariant probability subspace  $ \mathcal{U} \subset S(V) $
and for  any pairing $\mathcal{P} \in \Pi_{2l} $, let $ x( \mathcal{P})$ be a random Veronese tensor defined by a right hand side of (4). If representation of $G$ in $V$ is irreducible then  
\begin{equation}
		\mathbb{E}(  x( \mathcal{P}) ) = \frac{1}{n^l} I(\mathcal{P})    \nonumber
\end{equation}	 
\end{lemma}

\section{Expectations and Invariants}
From now on, unless explicitly stated otherwise, we assume that the group $G$ is a full orthogonal group $O(V) \approx O(n)$. 
For obvious reasons, everything that was said above remains valid in this special case. Note that $S(V) $ is an orbit of $O(V)$ and  hence the $O(V)$-invariant probability measure on a unit sphere $S(V)$ is uniquely inherited from the Haar measure on $
 O(n)$ (cf. e.g. \cite{Leng}).   
By direct computation (see  for example \cite{Folland}, \cite{Meckes}) one gets 
\begin{lemma}
	For   independent invariant random vectors $x,y$ defined on   $ S(V) $
	\begin{equation}
	\mu_{k,n} \; \equiv \; 	
	\mathbb{E}( <x,y>^{2k} ) = \frac{(2k-1)!!}{(n + 2k - 2)(n + 2k - 4) \cdots n} 
	\end{equation}
\end{lemma}
\noindent
 Let $m = 2k $. Recall, that 
the space $ \mathcal{I}_{m=2k} $ is spanned  
by the set of standard invariants  $ \mathcal{V} = \{  I(\mathcal{P}),   \mathcal{P} \in \Pi_m  \} $.  
Let 
\begin{equation}
	A_m =  \frac{1}{(m-1)!!} \sum_{\mathcal{P} \in \Pi_m}  I(\mathcal{P}_m)  \in \mathcal{I}_m   \nonumber
\end{equation}
be the average of the set of standard invariants in $\mathcal{I}_m$.
Denote the denominator on the right hand side of (5) by $ P(n,k ) $ and let $ \mathcal{G} $
denote the Gram matrix of the ordered standard set of invariants $\mathcal{V} $. 
\begin{theorem} ($m=2k$)
	
	\begin{enumerate}
		\item[(i)]  	The orthogonal projection of the Veronese surface $ \mathfrak{V}_m $ onto $ \mathcal{I}_m $ is just one point $ 	\mu_{k,n} A_m$. 
		\item[(ii)] Sum of the elements of any row of $\mathcal{G} $ is equal to $P(n,k)  $ and the average of elements of  $\mathcal{G} $  is equal to $ \frac{1}{	\mu_{k,n}}$  
		\item[(iii)]  Therefore, the  expectation of a random Veronese tensor $  x ^{\otimes,2k} $ is an average of elements of the standard basis of $O(V)$-invariants  divided by the average of the matrix elements of $\mathcal{G} $, in other words
		 	\begin{equation}
			\mathbb{E}(  x ^{\otimes,m=2k} ) =   \frac{1}{P(n,k)} \sum_{ \mathcal{P} \in \Pi_m} I(\mathcal{P}) \; = \; \frac{(2k-1)!!}{P(n,k)} A_m  \; \equiv \; \mu_{k,n} A_m\nonumber
		\end{equation}

\end{enumerate}

\end{theorem}
\noindent Proof.  Since $ \mathfrak{V}_m $ is an orbit of $S(V) $, its (equivariant) orthogonal projection onto $\mathcal{I}_m$ consists of just one point. To find this point, let    
\begin{equation}
	P = \sum_{ \mathcal{P} \in \Pi} \alpha_{ \mathcal{P}} I(\mathcal{P}_m)  \in \mathcal{I}_m   \nonumber
\end{equation}
be a projection of Veronese tensor $ a ^{\otimes,m}, \; a \in S(V) $ onto 
$\mathcal{I}_m$ where $  \alpha_{ \mathcal{P}} , \mathcal{P} \in \Pi_m$ are real numbers that must satisfy the "normal" equations
\begin{equation}
	<P -  a ^{\otimes,m}, \;  I(\mathcal{P}) > \; = \; 0 , \;  \text{for all } \mathcal{P} \in \Pi_m  \nonumber
\end{equation}
\noindent It is obvious that  $<a ^{\otimes,m},   I(\mathcal{P} )> \; = \; 1 $ for all $ \mathcal{P} \in \Pi_m $.  Using an appropriate order on the set of pairings   rewrite this system of  linear equations as 
\begin{equation}
	\mathcal{G}  \alpha = 1_{L}
\end{equation} 
where $L=(2k-1)!!, \;  \alpha$ is an unknown vector with coordinates 
$\alpha_{ \mathcal{P}}, \; \mathcal{P} \in \Pi_m $ and
$ 1_L$ is the $L$-dimensional vector of ones. Since $ \mathcal{V} $ is an $S_m$-orbit,  all row-sums  of $\mathcal{G}$ are equal to each other by Lemma 2. Denoting the unique row-sum of $ \mathcal{G}$ by $N$ we see that vector $(1/N) 1_L $ 
is a solution of the  system of equations (6). 
The  value of $N$ can be easily found by comparing (3) and (5).
Indeed, from Lemma 3 (with $m=2k$) and  (5), (6)  one gets
\begin{equation}
	L/N = \frac{1}{N^2}<\sum_{ \mathcal{P} \in \Pi_m} I(\mathcal{P}) , \sum_{ \mathcal{P} \in \Pi_m} I(\mathcal{P}) >  \; = \; L/P(n,k) \equiv \mu_{k,n}  \nonumber 
\end{equation}  
 It follows  that $ N= P(n,k)$ and that  the number $1/\mu_{k,n} $ is the average of elements of the Gram matrix $ \mathcal{G}$.  That proves (i) and (ii) while (iii) is equivalent to (i) by Lemma 1.
%
%
%
\newline

\noindent
\noindent It is now  a straightforward exercise to verify the following corollary
\begin{theorem}
 Let $\mathcal{P} = \{J_1, J_2, \cdots, J_l \} $ be a partition  of the set of indexes $ J= \{1,2, \cdots, m \} $,  let 
 	$ x_1, x_2, \cdots x_l \in S(V) $ be independent $O(V)$-invariant random vectors 
 and let $ X = x_1(J_1) \otimes \cdots \otimes x_l(J_l) $  be a generalized random Veronese tensor as in (4)	
 \begin{enumerate}
 	\item[(i)] If cardinality of at least one  of the sets $J_i, \; i = 1, \cdots l $ is odd then 
 	 $	\mathbb{E}( X) = 0 $
 \item[(ii)] If $ |J_i| = 2k_i, \; i = 1, \cdots l $,
 let $ \Pi(J_i) $ be the set of all $(2k_i - 1)!!$ pairings on the set $J_i, \; i = 1,\cdots l $ and let $\Pi = \Pi(J_1) \times \cdots \times \Pi(J_l) $.   Then in notation of Definition 1 (cf. Remark 1)
 \begin{equation}
 		\mathbb{E}( X) = \left( \prod_{i=1}^l \frac{(2k_i - 1)!!}{P(n,k_i)} \right)
 		\frac{1}{|\Pi|} \sum_{(\mathcal{P}_1, \;\cdots, \mathcal{P}_l)\in \Pi } 
 		I(\mathcal{P}_1 ) \otimes \cdots \otimes I(\mathcal{P}_l)   \nonumber
\end{equation}
\end{enumerate}
\end{theorem}
\begin{corollary} Let $  \mathcal{V}_0 = \{  I(\mathcal{P}_i),   \mathcal{P}_i \in \Pi' \subset  \Pi_m \}, \; i = 1,\cdots, m_0  \leq m  $ be an ordered basis of $ \mathcal{I}_{m=2k} $ and let $ \mathcal{G}_0$ be the corresponding Gram matrix of  $ \mathcal{V}_0 $. Let $ s_i, \; i = 1,\cdots , m_0 $  be the sum of elements of the $i$-th row of $ \mathcal{G}^{-1}_0$. Then the following identity holds
\begin{equation}
	\sum_{i=1}^{m_0} s_i I(\mathcal{P}_i) \; = \; \mu_{k,n} A_m \nonumber
\end{equation}
\end{corollary}

\noindent Proof.  Following the proof of Theorem  bear in mind that: (a) 
projection of the Veronese surface onto $ \mathcal{I}_m $ belongs to the  linear span of $\mathcal{V}_0$, and (b), 
 that $ (s_1, \cdots , s_{m_0})^T = \mathcal{G}^{-1}_0 1_{m_0} $
	
\begin{remark}
	Finding explicit values of $s_i$ in case of  $ n < k $ is probably a non-trivial problem. A method for selecting a basis from the standard set of invariants can be found in \cite{Alg} 
\end{remark}

\subsection{Gram Matrix of the Standard Set of Invariants  }
 For any two partitions $\mathcal{P}_1, \mathcal{P}_2 $
 denote their least upper bound by $  \mathcal{P}_1 \lor \mathcal{P}_2    $. The following useful statement is well known (cf. e.g. \cite{Weing_1} and \cite{gram}).
 
\begin{lemma}
	 Let 
	$\mathcal{P} = \{P_1, P_2, \cdots, P_k \} $ 
	and $\mathcal{Q} = \{Q_1, Q_2, \cdots, Q_k \} $ be two pairings   
	of the index set  $ J= \{1,2, \cdots, m=2k \} $, 
	let $\mathcal{R} =  \{R_1, R_2, \cdots, R_l \}  = \mathcal{P} \lor \mathcal{Q} $ and set $l = |\mathcal{P} \lor \mathcal{Q}|$. Then
	\begin{equation}
		<I(\mathcal{P}), I(\mathcal{Q})> \; = \; n^l
	\end{equation}

\end{lemma}
\noindent It is not hard to come up with a combinatorial proof of this statement. We prefer, however, to verify this fact by computing  expectations.
By Lemma 5, there are generalized random Veronese tensors 
$ x( \mathcal{P} ) = x_1(P_1) \otimes \cdots \otimes x_k(P_k) $ and $ y( \mathcal{Q}) = y_1(Q_1) \otimes \cdots \otimes y_k(Q_k  ) $ such that 
\begin{equation}
		\mathbb{E}(  x( \mathcal{P}) ) = \frac{1}{n^k} I(\mathcal{P}) \; \text{and}\;
			\mathbb{E}(  y( \mathcal{Q}) ) = \frac{1}{n^k} I(\mathcal{Q})	
\end{equation}
One can further assume that random unit vectors $x_i, y_i, i = 1, \cdots , k $ are pairwise independent. By separating  variables we have   
\begin{equation}
	\mathbb{E}( <x( \mathcal{P} ), y( \mathcal{Q}> ) = 
	\prod_{i=1}^{l} 	\int_{S(V)^{|\Theta_i|}} \prod_{x_a,y_b \in \Theta_i} <x_a,y_b> d\mu_i 
\end{equation}   
where  sets of variables $\Theta_i$ are pairwise disjoint, $ |\Theta_i| = |R_i|, i = 1, \cdots, l $ and every variable occurs exactly twice in each of the integrals over a product of $ |\Theta_i|$ unit spheres (denoted by $ S(V)^{|\Theta_i|} $).
 In vector coordinates, the $i$-th integral  splits into a sum of  integrals of monomials. Among these monomials there are exactly $n$ that are products of $|R_i|$ squares and the rest of the monomials will have at least one independent multiple of the form $x_i x_j,  i \neq j $ (two different coordinates of the same vector).
 It is easy to see (cf. Lemma 4 and \cite{irr}) that 
 square terms  will contribute a  
value of $ n^{-|R_i|}   $ each, while all other terms will vanish. 
 Hence, the right hand side of (9) evaluates to $n^{-(2k-l)} $ and then it follows   from (8)  that
\begin{equation}
\frac{1}{n^{2k\;-\;l}} = 		\mathbb{E}(  <x( \mathcal{P}), y( \mathcal{P}) >  ) = \frac{1}{n^{2k}}	<I(\mathcal{P}), I(\mathcal{Q})>  \nonumber
\end{equation}
which is the same as (7). 
\newline
\newline
\noindent  By Lemma 7  every entry of the Gram matrix $\mathcal{G} $ is a power of $ n $ and by  Theorem 1 (ii), the sum of entries of every row of $\mathcal{G} $ is equal to $P(n,k) = n^k + c_1 n^{k-1} + \cdots c_k n $ which is  a value of the polynomial of $n$ with fixed integer coefficients (that do not depend on $n$). Therefore,  every row of $\mathcal{G} $ has exactly one element equal to $n^k$, $c_1$ elements equal to $n^{k-1} $ and so on. It is easy to see, therefore, that the following statement is true 
\begin{corollary}
	For $ t = 1, \; 3, \;  5, \; \cdots , \; 2k -1 $ let $M_t$ be a $ t \times t $ matrix with all its diagonal elements equal to $n$ and all its off-diagonal elements equal to one. Let 
	$ M $ be the Kronecker product of matrices $ M_t, \; t = 1, \; 3, \; \cdots ,\; 2k -1  $.Then, the  rows of $ \mathcal{G} $ and $M $ are the same up-to a (row dependent) permutation. 
\end{corollary}
\noindent Another straightforward application of Lemmas 1, 5 and 7 yields  
\begin{corollary}
Let $G$ be a compact Lie group and let $V$ be a real $G$-representation.  
For any set of $2k$ independent random $G$-invariant vectors
\begin{equation}
x= \{x_1, x_2, \cdots,  x_k \}, \;  y = \{y_1, y_2, \cdots,  y_k \} \in \mathcal{U}  \nonumber
\end{equation}
defined on any  $G$-invariant probability subspace  $ \mathcal{U} \subset S(V) $
and for  any pairings $\mathcal{P}, \mathcal{Q} \in \Pi_{2k} $, let $ x( \mathcal{P}), y( \mathcal{P}) $ be corresponding  generalized random Veronese tensors in 
$ \mathcal{U}^{\otimes,m}$ defined as in   (4). If representation of $G$ in $V$ is irreducible then 
\begin{equation}
\mathbb{E}( <x( \mathcal{P} ),  y( \mathcal{Q}> ) = 	\frac{1}{n^{2k\;-\;| \mathcal{P}  \lor \mathcal{Q} |}}   \nonumber
\end{equation}
\end{corollary}
  
\subsection{  Weingarten Calculus  }
Let $X$ be a  random Haar-distributed matrix in $ O(V) \approx O(n) $. One of the standard problems addressed by
 Weingarten calculus (see e.g. \cite{Weing_1} and references therein)  is computation of moments  $	\mathbb{E}( x_{i_1,j_1} \cdots x_{i_m,j_m  } )$ where $x_{i,j} $  is an $ (i,j)$-entry of $X$. To apply Theorem 2 in this context, recall that 
 $ x_{i,j} = \; <Xe_{i}, e_{j} > $ and   hence (cf. e.g. \cite{Weing_1})   
\begin{equation}
 \mathcal{E} \equiv 	\mathbb{E}( x_{i_1,j_1} \cdots x_{i_m,j_m  } ) = 
 	\mathbb{E} (<Xe_{i_1}, e_{j_1} > \cdots <Xe_{i_m}, e_{j_m } ) >  
\end{equation} 
Let's $ \{i'_1, \cdots i'_l \} $ be the set of all  distinct first indices in (10). Suppose that $ i'_t $ occurs in the sequence 
$ i_1 \cdots i_m, \; r_t  $ times, $ \sum_{t=1}^l r_t = m $.
Let $\mathcal{P} = \{J_1, J_2, \cdots, J_l \} $ be a partition of the set of indexes $ J= \{1,2, \cdots, m \} $ such that $ J_1 $
consists of the first $ r_1 $ elements of $J, \; J_2 $ consists of the next $r_2$ elements of $ J $ and so on. Rearranging first indexes on the right hand side of (10) in accordance with the partition $\mathcal{P}$ gives rise to a  permutation of the corresponding sequence of second indexes in  dot-product terms. Let thus permuted second index sequence be $ \sigma = (j'_1, \cdots j'_m ) $. Split the sequence $\sigma$  into subsequences of $ \sigma_1, \sigma_2, \cdots, \sigma_l $ so that $\sigma_1$ contains first $r_1$ elements of $\sigma$, $\sigma_2$ contains next $r_2$ elements of $\sigma$ and so on. Further, 
for $ \sigma_i = ( a_1, a_2, \cdots, a_{r_i} ) $
set  $
	e(\sigma_i) = e_{a_{1}} \otimes e_{a_2} \otimes
		\cdots \otimes e_{a_{r_i}}
, \; i = 1, \cdots, l$.
    The next step is to form a random $O(n)$-invariant Veronese tensor (cf. (4))
\begin{equation}
	\mathcal{X}=  (Xe_{i'_1})(J_1) \otimes \cdots \otimes (Xe_{i'_l})(J_l)   
\end{equation} 
The  columns $Xe_{i'_t}, \; t = 1, \cdots l $ of the matrix $X$ are not independent as random vectors. However,  
these columns are pairwise orthogonal and hence 
$<I(\mathcal{Q}), \mathcal{X}> \; = \; 0 $ for any pairing $\mathcal{Q}$ on $J$ that 
is not a refinement of the partition  $\mathcal{P} $. Therefore, applying the
equivariant projector    
$ P^{O(V)} :  V^{\otimes,m} \rightarrow \mathcal{I}_m  $ from Lemma 1 to the generalized random Veronese tensor  (11) and arguing as in the proof of Theorem 1, we get 
\begin{equation}
	\mathbb{E}(\mathcal{X} ) = P^{O(V)}(\mathcal{X} ) =
	P_1^{O(V)}(  (Xe_{i'_1})(J_1)) \otimes \cdots \otimes 
	P_l^{O(V)}( (Xe_{i'_l})(J_l) ) 
\end{equation}
where $ 	P_t^{O(V)} $ is an equivariant projector 
$ P_t^{O(V)} :  V^{\otimes,k_t} \rightarrow \mathcal{I}_{k_t}  $,
and by Lemma 1 the right hand side of (12) is equal to 
\begin{equation}	
	\mathbb{E}( (Xe_{i'_1})(J_1))
	\otimes \cdots \otimes \mathbb{E}( (Xe_{i'_l})(J_l)) \nonumber
\end{equation}  
Assuming introduced notation we have a computational recipe.
\begin{theorem} (Cf. Theorem 2). Expectation $	\mathcal{E}$ in (10)  can be computed as follows:
\begin{enumerate}
		\item[(i)] If $|J_i| $ is odd for at least one of the sets $J_i, \; i = 1, \cdots l $ then 
		$	\mathcal{E} = 0 $
		\item[(ii)] If $ r_i = |J_i| = 2k_i, \; i = 1, \cdots l $,  
		let $ A_i $ be an average of standard $ 2k_i \;\;O(V)$-invariants defined by
		$(2k_i-1)!! $ pairings of the set $J_i, \; i = 1,\cdots l $. Then
		\begin{equation}
		\mathcal{E}  = \left( \prod_{i=1}^l \frac{(2k_i - 1)!!}{P(n,k_i)} \right) <A_1, e(\sigma_1)> < A_2, e(\sigma_2 )>\cdots < A_l, e(\sigma_l)>  \nonumber
		\end{equation}
	\end{enumerate}
\end{theorem} 
\begin{corollary} If indexes $ i_1, \cdots i_l $ are pairwise distinct then  
	 
	\begin{equation}
		\mathbb{E}( x_{i_1,j_1}^{2k_1} \cdots x_{i_l,j_l  }^{2k_l} ) = \prod_{i=1}^l \frac{(2k_i - 1)!!}{P(n,k_i)}   \nonumber
	\end{equation}
\end{corollary}

\end{document}